\documentclass{article}

\newtheorem{fed}{\textbf{Definition}}[section]
\newtheorem{thm}[fed]{\textbf{Theorem}}

\usepackage{amssymb,bbm,graphicx,epsfig,psfrag,epic,eepic,latexsym}
\usepackage{amsmath}
\usepackage{mathrsfs} 
\usepackage[dvips]{color}       

\begin{document}
\title{Existence of Bogomol'nyi solitons via Floer theory}
\author{Urs Frauenfelder \footnote{Guest at the Forschungsinstitut
f\"ur Mathematik at ETH Z\"urich (FIM)}}
\maketitle

\begin{abstract}
In this paper we use Floer homology to prove existence of 
solutions of a generalized Abelian Higgs equation introduced by
Schroers. 
\end{abstract}

\tableofcontents

\newpage

\section[Introduction]{Introduction}

A Floer homology is the Morse homology of an action functional
$\mathcal{A}$ defined on the infinite dimensional space
$C^\infty(N,M)$ of smooth maps from a manifold $N$ to a target manifold
$M$. The Floer chain complex is generated by critical points of 
$\mathcal{A}$ which are smooth maps from
$N$ to $M$ solving a PDE (respectively an ODE if $N$ is 1-dimensional).
The Floer boundary operator is defined by counting gradient flow lines
of $\mathcal{A}$. These are smooth maps from $N \times \mathbb{R}$
to $M$ again solving a PDE. In the special case where $N=S^1$ is
the circle and $\mathcal{A}$ is the classical Hamiltonian action
functional the critical points of $\mathcal{A}$ are the periodic orbits
of the Hamiltonian system. This was used by Floer in his proof of
the Arnold conjecture \cite{floer1,floer2,floer3,floer4}. 

Floer flow lines satisfy a Bogomol'nyi type selfduality equation 
and hence the solutions of various problems in physics can be
interpreted as Floer flow lines. In this paper we prove existence of
solutions of a generalized Abelian Higgs Equation discovered by 
Schroers \cite{schroers}. This equation is a special case
of the symplectic vortex equations introduced by Cieliebak, Gaio, Mundet,
and Salamon \cite{cieliebak-gaio-mundet-salamon, cieliebak-gaio-salamon}.
 
The method of our proof is the following. We consider an action functional
$\mathcal{A}$
which is invariant under the action of a finite group $\Gamma$
which acts freely on
the critical set of $\mathcal{A}$. This gives rise to two Floer homologies
where the chain complex of the first one is generated by critical points
of $\mathcal{A}$ and the chain complex of the second one is generated
by $\Gamma$-orbits of critical points of $\mathcal{A}$. We show
that the first Floer homology is trivial and the second one is non-trivial
which implies the existence of Floer flow lines.
\\ \\
This paper is organized as follows. In Section~\ref{ffl} we
show how Floer flow lines can be interpreted as selfduality 
equations and how they are related to problems in physics.
In Section~\ref{statement} we introduce our Floer action functional
and state our main theorem. In Section~\ref{proof} we first explain
how one can find Floer flow lines if one is able to compute two
Floer homologies. We then introduce our two Floer homologies and 
compute them. We finally consider the classical vortex equations
as an example.
\\ \\
\textbf{Acknowledgements: }I would like to thank Dietmar Salamon
for useful discussions. This paper was written during my stay at
FIM of ETH Z\"urich. I wish to thank FIM for its kind
hospitality.

\section[Floer flow lines as selfduality equations]{Floer flow lines
as selfduality equations}\label{ffl}

Assume that $(M,g)$ is a Riemannian manifold and $f \in C^\infty(M)$ is a
Morse function on $M$. For two critical points $c_-, c_+ \in\mathrm{crit}(f)$
let  
$$\Omega(c_-,c_+)=\{x \in C^\infty(\mathbb{R},M):
\lim_{s \to \pm \infty}x(s)=c_\pm\}$$
be the space of paths from $c_-$ to $c_+$. The energy functional 
$E \colon \Omega(c_-,c_+) \to \mathbb{R} \cup \{\infty\}$ is defined by
$$E(x)=\frac{1}{2}
\int_{\mathbb{R}}\Big(||\partial_s x||_g^2+||\nabla_g f(x)||_g^2\Big)ds.
$$
The energy functional satisfies the following Bogomol'nyi equation
\cite{bogomolnyi},
\begin{eqnarray}\label{bogomolnyi}
E(x)&=&\frac{1}{2}\int_{\mathbb{R}}||\partial_s x+\nabla_g f(x)||_g^2 ds
+f(c_-)-f(c_+)\\ \nonumber
&=&\frac{1}{2}\int_{\mathbb{R}}||\partial_s x-\nabla_g f(x)||_g^2 ds
+f(c_+)-f(c_-).
\end{eqnarray}
It follows from (\ref{bogomolnyi}) that downwards and upwards gradient flow
lines, i.e.\ solutions of the following first order ODE 
\begin{equation}\label{dugrad}
\partial_s x=\pm \nabla_g f(x)
\end{equation}
are the absolute minimizers of the energy functional. Note that the
Euler-Lagrange equations for $E$ are given by the following second order
ODE
\begin{equation}\label{euler-l}
\nabla^2_s x=H_f(x)\nabla_g f(x)
\end{equation} 
where $H_f$ denotes the Hessian of $f$. Of course solutions of
(\ref{dugrad}) are also solutions of (\ref{euler-l}). This can be
directly checked by differentiating (\ref{dugrad}).
Equations describing the absolute minima of a functional which 
are first order
while the Euler-Lagrange equations of the functional
are second order are called
\textbf{(anti)-selfduality equations}. 

Selfduality equations are of great importance in physics. We 
show how some of the selfduality equations physicists are interested
in can be interpreted as gradient flow lines of an action functional
on an infinite dimensional space. The basic example are the
classical vortex equations on the cylinder, see \cite{jaffe-taubes}
as a basic reference. To describe them
recall that the standard $S^1$-action on the complex plane $\mathbb{C}$
given by
$$z \mapsto e^{i\theta}z, \quad e^{i\theta} \in S^1$$
is Hamiltonian with respect to the standard symplectic structure
$\omega=dx \wedge dy$ on $\mathbb{C}$. A moment map for 
this action is given by
$$\mu(z)=-\frac{i}{2}|z|^2+\frac{i}{2} \in i\mathbb{R}=\mathrm{Lie}(S^1).$$
The vortex equations 
on the cylinder for a triple $(v,\eta,\zeta) \in C^\infty(\mathbb{R}
\times S^1, \mathbb{C}\times i\mathbb{R}\times i\mathbb{R})$
are
\begin{equation}\label{ve}
\left.\begin{array}{c}
\partial_s v+\zeta v+i(\partial_t v+\eta v)=0,\\
\partial_s \eta-\partial_t \zeta+\mu(v)=0.
\end{array}\right\}
\end{equation}
If one thinks of $A=\zeta ds+\eta dt$ as a connection on 
the trivial bundle over the cylinder then
the second equation in (\ref{ve}) 
reads
$$*F_A+\mu(v)=0$$
where $*$ denotes the Hodge operator on the cylinder $\mathbb{R}\times S^1$
and $F_A$ is the curvature of the connection $A$. The gauge group
$$\mathcal{G}=C^\infty(\mathbb{R}\times S^1, S^1)$$
acts on the solutions of (\ref{ve}) via
$$g_*(v,\eta,\zeta)=(gv,\eta-g^{-1}\partial_t g, \zeta-g^{-1}\partial_s g).$$
In particular, every solution of (\ref{ve}) is gauge equivalent
to a solution of (\ref{ve}) with $\zeta\equiv 0$, i.e.\ a pair
$(v,\eta) \in C^\infty(\mathbb{R}\times S^1, \mathbb{C}\times i\mathbb{R})$
satisfying
\begin{equation}\label{vet}
\left.\begin{array}{c}
\partial_s v+i(\partial_t v+\eta v)=0,\\
\partial_s \eta+\mu(v)=0.
\end{array}\right\}
\end{equation}
The equations (\ref{vet}) are called vortex equations in temporal gauge.
The gauge group 
$$\mathcal{H}=C^\infty(S^1,S^1)$$
acts on solutions of (\ref{vet}) by
$$h_*(v,\eta)=(hv,\eta-h^{-1}\partial_t h)$$
and the map $(v,\eta) \mapsto (v,\eta,0)$ from solutions of
(\ref{vet}) to solutions of (\ref{ve}) induces a bijection
$$\{\mathrm{solutions\,\,of \,\,(\ref{vet})}\}/\mathcal{H}
\cong\{\mathrm{solutions\,\,of\,\,(\ref{ve})}\}/\mathcal{G}.$$
The equation (\ref{vet}) can be interpreted as the gradient flow
equation of an action functional $\mathcal{A}$ defined on
an infinite dimensional space $\mathscr{L}$. Set
$$\mathscr{L}=C^\infty(S^1,\mathbb{C}\times i\mathbb{R})$$
and define $\mathcal{A}\colon \mathscr{L} \to \mathbb{R}$ as
$$\mathcal{A}(v,\eta)=\int_0^1\lambda(v)\partial_t v
+\int_0^1\langle \mu(v),\eta \rangle dt$$
where $\lambda=y dx$ is the Liouville one-form satisfying
$d\lambda=-\omega$. Then the gradient flow lines of
$\mathcal{A}$ with respect to the $L^2$-metric $g_{L^2}$ on
$\mathscr{L}$, i.e.\
$$\partial_s(v,\eta)+\nabla_{g_{L^2}}\mathcal{A}(v,\eta)=0$$
are precisely the equations (\ref{vet}).
\\
Another example are the self-dual Chern-Simons vortices discovered
by Hong-Kim-Pac and Jackiw-Weinberg, see 
\cite{hong-kim-pac, jackiw-weinberg}. They read
\begin{equation}\label{csv}
\left.\begin{array}{c}
\partial_s v+\zeta v+i(\partial_t v+\eta v)=0,\\
\partial_s \eta-\partial_t \zeta+|v|^2\mu(v)=0.
\end{array}\right\}
\end{equation}
In temporal gauge they are given by
\begin{equation}\label{csvt}
\left.\begin{array}{c}
\partial_s v+i(\partial_t v+\eta v)=0,\\
\partial_s \eta+|v|^2\mu(v)=0.
\end{array}\right\}
\end{equation}
The equations (\ref{csvt}) are again gradient flow lines of
$\mathcal{A}$, but with respect to a different metric, namely the
following warped product metric on 
$\mathscr{L}=C^\infty(S^1,\mathbb{C})\times C^\infty(S^1,\times i\mathbb{R})$
given for $(v,\eta) \in \mathscr{L}$ and 
$(\hat{v}_1,\hat{\eta}_1), (\hat{v}_2,\hat{\eta}_2) \in 
T_{(v,\eta)}\mathscr{L} \cong \mathscr{L}$ by the formula
$$g(v,\eta)\big((\hat{v}_1,\hat{\eta}_1),(\hat{v}_2,\hat{\eta}_2)\big)=
\int_0^1\langle \hat{v}_1,\hat{v}_2 \rangle dt+
\int_0^1 \frac{1}{|v|^2}\langle \hat{\eta}_1,\hat{\eta}_2 \rangle dt.$$

\section[Statement of the main result]{Statement of the main result}
\label{statement}

Assume that for $k\leq n$ the torus 
$T^k=\{e^{iv}: v \in \mathbb{R}^k\}$ 
acts on the complex vector space
$\mathbb{C}^n$ via the action
$$\rho(e^{iv})z=e^{iAv}z, \quad z \in \mathbb{C}^n, \,\,v \in \mathbb{R}^k$$
for some $(n\times k)$-matrix $A$ with 
integer entries. We endow the 
Lie algebra of the torus
$$\mathrm{Lie}(T^k)=\mathfrak{t}^k=i\mathbb{R}^k$$
with its standard inner product. 
The action of the torus
on $\mathbb{C}^n$ is Hamiltonian with respect to the standard symplectic
structure $\omega=\sum_{i=1}^n dx_i \wedge dy_i$.
Let $\tau$ be an element of the Lie algebra $\mathfrak{t}^k$.  
Denoting by $A^T$ the transposed matrix of $A$ a moment map 
$\mu \colon \mathbb{C}^n \to \mathfrak{t}^k$ is given by
\begin{equation}\label{moment1}
\mu(z)=-i A^Tw-\tau, \quad w=\frac{1}{2}\left(\begin{array}{c}
|z_1|^2 \\
\vdots\\
|z_n|^2
\end{array}\right),
\end{equation}
i.e.
$$d\langle \mu, \xi \rangle=\iota_{X_\xi}\omega, \quad \xi \in 
\mathfrak{t}^k$$
for the vector field $X_\xi$ on $\mathbb{C}^n$ given by the
infinitesimal action
$$X_\xi(z)=\dot{\rho}(\xi)(z), \quad z \in \mathbb{C}^n.$$
We will assume throughout
this paper the following hypothesis, 
\begin{description}
 \item[(H)] The moment map $\mu$ is proper, $\mu^{-1}(0)$ is
 not empty, and $T^k$ acts freely on it.
\end{description}
It follows from (H) that the Marsden-Weinstein quotient
$$\mathbb{C}^n//T^k=\mu^{-1}(0)/T^k$$
is a compact symplectic manifold of dimension
$$\mathrm{dim}(\mathbb{C}^n//T^k)=2(n-k),$$
where the symplectic structure is
induced from the standard symplectic structure on $\mathbb{C}^n$.

Let $H_t \in C^\infty(\mathbb{C}^n)$ for $t \in [0,1]$ be a smooth
family of $T^k$-invariant Hamiltonians. We denote by
$$dH_t=-\omega(X_{H_t},\cdot)$$
the Hamiltonian vector field of $H_t$.
Furthermore, let $J_t$ for $t \in [0,1]$
be a smooth family of $T^k$-invariant $\omega$-compatible almost
complex structures, i.e.
$$g_t(\cdot,\cdot)=\omega(\cdot,J_t\cdot)$$
is a Riemannian metric  on $\mathbb{C}^n$.
In this paper we are proving existence of
solutions $(v,\eta) \in C^\infty(\mathbb{R}\times [0,1],\mathbb{C}^n\times
\mathfrak{t}^k)$ of the following PDE
\begin{equation}\label{sve}
\left.\begin{array}{c}
\partial_s v+J_t(v)(\partial_t v+X_\eta(v)-X_{H_t}(v))=0,\\ 
\partial_s \eta+\mu(v)=\tau,\\
v(0), v(1) \in \mathbb{R}.
\end{array}\right\}
\end{equation}
Schroers \cite{schroers} discovered these equations for 
the special case where $J_t$ is the standard complex structure
given by multiplication with $i$ and the Hamiltonian $H$ vanishes
identically, see also \cite{yang}. 
Note that for the standard circle action on
1-dimensional complex space they reduce to the classical vortex
equations. 
They are examples of the symplectic vortex equations in temporal gauge
on the cylinder. The symplectic vortex equations
can be defined more generally for Hamiltonian group actions
on a symplectic manifold. 
They were introduced by
Cieliebak, Gaio, Mundet, and Salamon in 
\cite{cieliebak-gaio-mundet-salamon, cieliebak-gaio-salamon}.

Abbreviate
$$T^k_{\mathbb{R}^n}=\{g \in T^k: \rho(g)\mathbb{R}^n=\mathbb{R}^n\}$$
the isotropy subgroup of the Lagrangian submanifold $\mathbb{R}^n$.
Since by our standing assumption (H) there is a point
$x \in \mathbb{R}^n$ which is fixed only by the identity in $T^k$ 
it follows that
$$T^k_{\mathbb{R}^n} \cong \mathbb{Z}_2^k.$$
The gauge group 
$$\mathcal{H}=\{h \in C^\infty([0,1],T^k): h(0), h(1) \in T^k_{\mathbb{R}^n}\}
$$
acts on the solutions of (\ref{sve}) by
$$h_*(v,\eta)=(\rho(h)v, \eta-h^{-1}\partial_t h).$$  
We define the energy of a solution of (\ref{sve}) by
$$E(v,\eta)=\int_{\mathbb{R} \times S^1}\big(g_t(\partial_s v,\partial_s v)
+\langle \partial_s \eta,\partial_s \eta \rangle\big)ds dt.$$
We finally introduce the space of admissible families of
Hamiltonians $\mathscr{H}$ 
and the space of admissible families of almost complex structures
$\mathscr{J}$. The space $\mathscr{H}$ consists 
of smooth families of
$T^k$ invariant Hamiltonians $H_t$ for which there exists
a compact set $K=K(H) \subset \mathbb{C}^n$ such that
the support of $H_t$ for every $t \in [0,1]$ is contained in
$K$. The space $\mathscr{J}$ consists of smooth families
of $T^k$-invariant $\omega$-compatible almost complex structures
$J_t$ such that there exists a compact set $K=K(J) \subset \mathbb{C}^n$
such that $J_t$ equals outside of $K$ 
the standard complex structure, i.e.\ the
multiplication by $i$.\\
Our main result is the following theorem.
\\ \\
\textbf{Theorem A:} 
\emph{For generic pairs $(H,J) \in \mathscr{H}\times \mathscr{J}$ 
there exists at least two not gauge equivalent solutions 
of (\ref{sve}) whose energy is positive and finite.}

\section[Proof of Theorem A]{Proof of Theorem A}\label{proof}

\subsection[A general principle for finding Floer flow lines]{
A general principle for finding Floer flow lines}

Let $\mathscr{L}$ be a space on which a finite group $\Gamma$ acts and
let $\mathcal{A} \colon \mathscr{L} \to \mathbb{R}$ be a $\Gamma$-invariant
action functional all whose critical points are nondegenerate. 
Assume that the action of $\Gamma$ on
$\mathrm{crit}(\mathcal{A})$ is free. We then can define two Floer chain
complexes $CF(\mathscr{L},\mathcal{A})$ and
$CF(\mathscr{L}/\Gamma, \mathcal{A})$ where the first one 
is generated by critical points of $\mathcal{A}$ and
the second one is generated by $\Gamma$-orbits of critical points
of $\mathcal{A}$. Obviously,
$$\mathrm{rk} \big(CF(\mathscr{L},\mathcal{A})\big)=|\Gamma| \cdot
\mathrm{rk} \big(CF(\mathscr{L}/\Gamma,\mathcal{A})\big).$$
Suppose that we are able to define Floer homologies 
$HF(\mathscr{L},\mathcal{A})$ and $HF(\mathscr{L}/\Gamma,\mathcal{A})$
out of the chain complexes above. We are now in position to formulate
our principle.
\\ \\
\textbf{Principle: } \emph{Assume that
$\mathrm{rk} \big(HF(\mathscr{L},\mathcal{A})\big)\neq |\Gamma| \cdot
\mathrm{rk} \big(HF(\mathscr{L}/\Gamma,\mathcal{A})\big).$
Then there exists a Floer flow line of $\mathcal{A}$
of positive, finite energy.}
\\ \\
To illustrate our principle we consider an example were
$\mathscr{L}$ is a finite dimensional manifold. 
Let $\mathscr{L}$ be the two-sphere $S^2$ and let $\Gamma=\mathbb{Z}_2$
act on $S^2$ as the antipodal involution. Then 
the quotient $S^2/\Gamma$
is $\mathbb{RP}^2$. We consider the Morse homology with
$\mathbb{Z}_2$-coefficients.
In this case the assumption of our principle is
satisfied, i.e.\
$$\mathrm{rk}\big(HM(S^2;\mathbb{Z}_2)\big) \neq 
|\Gamma|\cdot \mathrm{rk}\big(HM(\mathbb{RP}^2;\mathbb{Z}_2)\big).$$
To understand what is going on we consider the
standard Morse function on $\mathbb{RP}^2$ which has 
one maximum, one critical point of index one, and one minimum,
two Morse flow lines from the maximum to the critical point of index one as
well as two Morse flow lines from the critical point of index one to the 
minimum. Since we are considering Morse homology with
$\mathbb{Z}_2$-coefficients the boundary operator on the Morse chain
complex on $\mathbb{RP}^2$ is zero. However, on $S^2$ the
two flow lines split off and the boundary operator of the Morse complex
on $S^2$ is not trivial anymore.  

\medskip
\begin{center}\begin{picture}(0,0)%
\includegraphics{s2.pstex}%
\end{picture}%
\setlength{\unitlength}{3947sp}%
\begingroup\makeatletter\ifx\SetFigFont\undefined%
\gdef\SetFigFont#1#2#3#4#5{%
  \reset@font\fontsize{#1}{#2pt}%
  \fontfamily{#3}\fontseries{#4}\fontshape{#5}%
  \selectfont}%
\fi\endgroup%
\begin{picture}(5513,3662)(451,-4540)
\end{picture}
\end{center}
\medskip

\subsection[Two Floer homologies]{Two Floer homologies}

We first show how solution of (\ref{sve}) arise as flow lines
of an action functional.
We introduce the following path space
$$\mathscr{P}=\{(v,\eta) \in C^\infty([0,1],\mathbb{C}^n \times 
\mathfrak{t}^k) \colon v(0), v(1) \in \mathbb{R}^n\}.$$
For a smooth family of $T^k$-invariant Hamiltonians $H_t$
we define the action functional $\mathcal{A}_H \colon
\mathscr{P} \to \mathbb{R}$ by
$$\mathcal{A}_H(v,\eta)=\int_0^1 \lambda(v)\partial_t v
-\int_0^1 H_t(v(t))dt +\int_0^1\langle \mu(v),\eta\rangle dt$$
where $\lambda=\sum_{i=1}^n y_i dx_i$ is the 
Liouville 1-form satisfying $d \lambda=-\omega$.
Recall the gauge group 
$$\mathcal{H}=\{h \in C^\infty([0,1],T^k): h(0), h(1) \in T^k_{\mathbb{R}^n}\}.
$$
It acts on $\mathscr{P}$ by
$$h_*(v,\eta)=(\rho(h)v, \eta-h^{-1}\partial_t h).$$
Note that the differential of $\mathcal{A}$ is invariant under the
action of $\mathcal{H}$ on $\mathscr{P}$. For
a smooth family of $T^k$-invariant $\omega$-compatible almost
complex structures $J_t$ define the metric $g_J$ on
$\mathscr{P}$ for $(v,\eta) \in \mathscr{P}$ and
$(\hat{v}_1,\hat{\eta}_1), (\hat{v}_2,\hat{\eta}_2) \in
T_{(v,\eta)}\mathscr{P} \cong \mathscr{P}$ by
$$g_J(v,\eta)\big((\hat{v}_1,\hat{\eta}_1),(\hat{v}_2,\hat{\eta}_2)\big)
=\int_0^1\omega(\hat{v}_1, J_t(v)\hat{v}_2)dt 
+\int_0^1\langle \hat{\eta}_1,\hat{\eta}_2\rangle dt.$$
Then the gradient flow lines of $\mathcal{A}_H$ with respect to
$g_J$ are precisely the solutions of (\ref{sve}).

Recall that $\mathscr{H}$ denotes the space of families of admissible
Hamiltonians, i.e.\ torus invariant Hamiltonians with compact support,
and $\mathscr{J}$ denotes the space of families of admissible
almost complex structures, i.e\ families of torus invariant $\omega$-compatible
almost complex structures which agree outside of a compact set with
the standard complex structure on $\mathbb{C}^n$. For
generic $(H,J) \in \mathscr{H} \times \mathscr{J}$ it
was shown in \cite{frauenfelder} that the Floer homologies
$HF(\mathscr{P},\mathcal{A}_H,g_J;\mathbb{Z}_2)$ and 
$HF(\mathscr{P}/T^k_{\mathbb{R}^n},\mathcal{A}_H,g_j;\mathbb{Z}_2)$ 
are well defined. Theorem A follows from the following theorem.

\begin{thm} $HF(\mathscr{P},\mathcal{A}_H,g_J;\mathbb{Z}_2)=0$ and
$HF(\mathscr{P}/T^k_{\mathbb{R}^n},\mathcal{A}_H, g_J;\mathbb{Z}_2) 
\neq 0$.
\end{thm}
\textbf{Proof: }That 
$HF(\mathscr{P}/T^k_{\mathbb{R}^n},\mathcal{A}_H,g_J;\mathbb{Z}_2)$
is not zero follows from the results in \cite{frauenfelder}.
There it was shown that the homology
$HF(\mathscr{P}/T^k_{\mathbb{R}^n},\mathcal{A}_H,g_J;\mathbb{Z}_2)$
is given by the Morse homology with coefficients in $\mathbb{Z}_2$
of the induced Lagrangian in the Marsden-Weinstein quotient
$$\bar{L}=T^k(\mu^{-1}(0) \cap \mathbb{R}^n)/T^k
\cong (\mu^{-1}(0)\cap \mathbb{R}^n)/T^k_{\mathbb{R}^n}$$
tensored with a Novikov ring which is itself an infinite dimensional
vector space over the field $\mathbb{Z}_2$. Since $\mu^{-1}(0)$
is not empty it follows that $\bar{L}$ is not empty and
hence $HF(\mathscr{P}/T^k_{\mathbb{R}^n},\mathcal{A}_H,g_J;\mathbb{Z}_2)$
is not zero. 

To show that $HF(\mathscr{P},\mathcal{A}_J,g_J;\mathbb{Z}_2)$
is zero we make use of the fact that there exists a Hamiltonian
isotopy $\phi_H$ of $\mathbb{C}^n$ generated by a smooth family
of Hamiltonians $H_t$ with support in a compact set
$K \subset \mathbb{C}^n$ such that
$$\phi_H(\mathbb{R}^n) \cap \mu^{-1}(0)=\emptyset.$$
To define the corresponding Floer homotopy care has to be taken since
the Hamiltonians $H_t$ are not $T^k$-invariant anymore. 
We therefore define the Floer homotopy not on $\mathscr{P}$ but
on a Coulomb section $\mathscr{P}_c \subset \mathscr{P}$.

Here is how this works. There is a natural splitting
$$\mathcal{H}\cong \mathcal{H}_0 \times \mathbb{Z}^k \times 
T^k_{\mathbb{R}^n}$$
where the contractible infinite dimensional Lie group
$\mathcal{H}_0$ is given by elements $h=\exp(\xi) \in \mathcal{H}$
satisfying $\xi(0)=\xi(1)=0 \in \mathfrak{t}^k$. 
The gauge group $\mathcal{H}_0$ acts freely on $\mathscr{P}$.
Moreover, for each path $\eta \in C^\infty([0,1],\mathfrak{t}^k)$
there is a unique element $h_\eta \in \mathfrak{t}^k$ which puts
$\eta$ into Coulomb gauge
$$d^*\big((h_\eta)_*\eta\big)=-\partial_t\big((h_\eta)_*\eta)=0,$$
i.e. $(h_\eta)_*\eta$ is constant. Hence we may think of
$$\mathscr{P}_c=\{v \in C^\infty([0,1],\mathbb{C}^n)\colon v(0),v(1)
\in \mathbb{R}^n\}\times \mathfrak{t}^k$$ 
as a global section in the principal $\mathcal{H}_0$-bundle
$\mathscr{P}$, i.e. the following diagram is commutative
\\ \\
\setlength{\unitlength}{1.5cm}
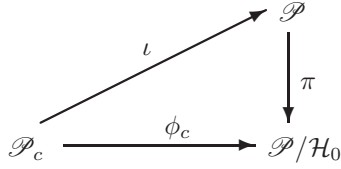
\begin{picture}(5,2)\thicklines
 \put(2,0){$\mathscr{P}_c$}
 \put(2.3,0.3){\vector(2,1){2}}
 \put(3.2,0.9){$\iota$}
 \put(2.5,0.1){\vector(1,0){1.7}}
 \put(3.4,0,2){$\phi_c$}
 \put(4.3,0){$\mathscr{P}/  \mathcal{H}_0$}
 \put(4.4,1.2){$\mathscr{P}$}
 \put(4.5,1.1){\vector(0,-1){0.8}}
 \put(4.6,0.6){$\pi$}
\end{picture}
\\ \\
where $\iota$ is the canonical inclusion and $\phi_c$ is the
isomorphism given by Coulomb gauge. 

The $L^2$-metric $g_{L^2}$ on $\mathscr{P}$ is $\mathcal{H}_0$-invariant, hence
it induces a quotient metric $[g_{L^2}]$ on $\mathscr{P}/\mathcal{H}_0$
and we denote by
$$g_c=\phi_c^*[g_{L^2}]$$
its pullback on $\mathscr{P}_c$. $\mathcal{H}_0$-equivalence classes
of gradient flow lines of $\mathcal{A}_H$ with respect to
$g_{L^2}$ are in natural one-to-one correspondence with gradient 
flow lines of the restriction of $\mathcal{A}_H$ to
$\mathscr{P}_c$ with respect to the metric $g_c$. 
To compute the metric we introduce the following notation. 
For $z \in \mathbb{C}^n$ 
denote by $L_z \colon \mathfrak{t}^k \to
T_z \mathbb{C}^n \cong \mathbb{C}^n$ the linear map
defined by
$$L_z \eta=X_\eta(z)$$
and by $L_z^*$ its dual with respect to the inner products on
$\mathfrak{t}^k$ and on $\mathbb{C}^n$. Then for
$(\hat{v}_1,\hat{\eta}_1), (\hat{v}_2,\hat{\eta}_2) \in
T_{(v,\eta)}\mathscr{P}_c$ where 
$(v,\eta) \in \mathscr{P}_c$ the metric $g_c$ is given by
\begin{eqnarray*}
g_c(v,\eta)\big((\hat{v}_1,\hat{\eta}_1),(\hat{v}_2,\hat{\eta}_2)\big)
&=&\int_0^1\langle \hat{v}_1-L_v\xi_{v,\hat{v}_1},\hat{v}_2 \rangle dt
+\langle \hat{\eta}_1,\hat{\eta}_2\rangle\\
&=&\int_0^1 \langle \hat{v}_1,\hat{v}_2-L_v\xi_{v,\hat{v}_2}\rangle dt+
\langle \hat{\eta}_1,\hat{\eta}_2 \rangle
\end{eqnarray*}
where $\xi_{v,\hat{v}}$ is determined by
$$L_v^*\hat{v}-L_v^*L_v \xi_{v,\hat{v}}+\partial^2_t
\xi_{v,\hat{v}}=0,\quad \xi_{v,\hat{v}}(0)=
\xi_{v,\hat{v}}(1)=0.$$

By abuse of notation we will denote by $\mathcal{A}_H$ also the
restriction of $\mathcal{A}_H$ to the Coulomb section $\mathscr{P}_c$.
In the following we will drop the assumption that
$H_t$ is $T^k$-invariant. To compute the gradient of 
$\mathcal{A}_H$ with respect to $g_c$ it is useful to
abbreviate for $v \in C^\infty([0,1],\mathbb{C}^n)$
\begin{eqnarray*}
\bar{\mu}(v)&=&\int_0^1 \mu(v(t))dt,\\
\kappa_{H}(v)(t)&=&\int_0^t L^*_{v(\tau)} \nabla H_\tau(v(\tau)) d\tau,\\
\bar{\kappa}_H(v)&=&\int_0^1 \kappa_H(v)(t)dt.
\end{eqnarray*}
Note that $\kappa_H$ and $\bar{\kappa}_H$ vanish
identically if $H$ is $T^k$-invariant. 
Using this notation we get for the gradient
$$\nabla_{g_c}\mathcal{A}_H(v,\eta)
=
\left(\begin{array}{c}
L_v \xi_v+i(\partial_t v+L_v\eta-X_{H_t}(v))\\
\bar{\mu}(v)
\end{array}\right)
$$
where $\xi_v$ is given by
$$\xi_v(t)=\int_0^t\big(\mu(v(\tau))+\kappa_{H}(v)(\tau)\big)d\tau
-t\big(\bar{\mu}(v)+\bar{\kappa}_H(v)\big).$$
Hence gradient flow lines satisfy the following equation
\begin{equation}
\left.\begin{array}{c}
\partial_s v +\dot{\rho}(\xi_v)+i(\partial_t v+\dot{\rho}(\eta))=0,\\
\partial_s \eta-\partial_t \xi_v+\mu(v)+\kappa_H(v)-\bar{\kappa}_H(v)=0.
\end{array}\right\}
\end{equation}
Up to the factor $\kappa_H(v)-\bar{\kappa}_H(v)$ which measures how far
is $H$ from being $T^k$-invariant these equations are examples of
symplectic vortex equations. However, since $H$ has compact support
there exists a constant $c=c(H)<\infty$ not depending on $v$ such that
$$||\kappa_H(v)-\bar{\kappa}_H(v)|| \leq c(H).$$
So the compactness result in \cite{cieliebak-gaio-mundet-salamon}
can be adapted to our situation to prove that the Floer homotopy is
well defined.  In our situation the proof simplifies considerably, though.
This is firstly due to the fact that our gauge group is abelian so
that we do not need the full strength of Uhlenbeck's compactness
theorem. Secondly thanks to the now established abstract perturbation
techniques we do not have to bother about transversality any more 
and can work with the standard complex structure on $\mathbb{C}^n$.
Since there are no derivatives of the standard complex structure
the bubbling analysis can be avoided by directly using
the elliptic estimate for the $\bar{\partial}$-operator. 
This proves the theorem. \hfill $\square$

\subsection[An example: The vortex equations]{An example: The vortex
equations} 

We consider the special case of the standard $S^1$-action on
$\mathbb{C}$. We set the Hamiltonian equal to zero and consider
the standard complex structure on $\mathbb{C}$.
Then the gradient flow lines are solutions of the classical
vortex equations on the strip with Lagrangian boundary condition.
The connected components of the critical set of $\mathcal{A}$
are parametrised by $\mathbb{Z} \times \mathbb{Z}_2$. It was
shown by Jaffe and Taubes \cite{jaffe-taubes} that
between critical components whose difference in the
$\mathbb{Z}$-factor is one there is exactly one solution
up to gauge equivalence. 

\medskip
\begin{center}\begin{picture}(0,0)%
\includegraphics{s1.pstex}%
\end{picture}%
\setlength{\unitlength}{3947sp}%
\begingroup\makeatletter\ifx\SetFigFont\undefined%
\gdef\SetFigFont#1#2#3#4#5{%
  \reset@font\fontsize{#1}{#2pt}%
  \fontfamily{#3}\fontseries{#4}\fontshape{#5}%
  \selectfont}%
\fi\endgroup%
\begin{picture}(5567,5316)(451,-5965)
\end{picture}
\end{center}
\medskip

Hence if one considers the Floer complex modulo
$\mathbb{Z}_2$, the
isotropy group of the Lagrangian $\mathbb{R}$, then
the Floer boundary operator vanishes and the homology is
nontrivial. On the other if one does not mod out by $\mathbb{Z}_2$,
then the same phenomenon as in the case of the standard
Morse function on $\mathbb{RP}^2$ happens. The flow lines split off
and the Floer homology is zero.

\end{document}